\documentclass[12pt]{article}

\usepackage{amssymb, amsmath, amscd, amsfonts, amsthm}
\usepackage[dvips,hypertex,colorlinks]{hyperref}%,backref,bookmarksnumbere
\input{epsf}

\begin{document}

\vspace*{2cm}
\begin{center}
{\bf\large INEQUALITIES FOR THE POLAR DERIVATIVE OF A POLYNOMIAL}\\[1cm]
{\bf M. Shakeri, M. Bidkham and M.Eshaghi Gordji}\\[0.2cm]
Department of Mathematics, Faculty of Natural Sciences \\
Semnan University, Semnan, Iran. \\
{\tt E-mail: mbidkham@semnan.ac.ir, madjid.eshaghi@gmail.com} \\[0.5cm]

\end{center}

\begin{abstract}
In this paper we obtain new results concerning maximum modulus of
the polar derivative of a polynomial with restricted zeros. Our
results generalize and refine upon the results of Aziz and Shah [An
integral mean estimate for polynomial, Indian J. Pure Appl. Math. 28
(1997) 1413--1419] and Gardner, Govil and Weems [Some result
concerning rate of growth of polynomials, East J. Apporox. 10(2004)
301--312].
\end{abstract}

\section*{\small 1. INTRODUCTION AND STATEMENT OF RESULTS}

The problems in the analytic theory of polynomials concerning
derivative of the polynomials have been frequently investigated.
Over many decades, a large number of research papers, e.g,
\cite{A1,A2, A3, A4, B, G2,M} have been published.

If $p(z)=\sum_{m=0}^{n}a_mz^m$ is a polynomial of degree$n$, then
\begin{equation}
\max_{|z|=1} |p'(z)|\leq n \max_{|z|=1} |p(z)|.
\end{equation}
The above inequality, which is an immediate consequence of
Bernstein's inequality on the derivative of a trigonometric
polynomial is best possible with equality holding for the polynomial
$p(z) = \lambda z^n$, $\lambda$ being a complex number.

\vrule height .001cm depth .002cm width 4.5cm\\
%{\bf ---------------------------------}\\
{\scriptsize 2000  Mathematics Subject Classification. 30A10, 30C15,
26C10, 12D10, 30D15.
%\vspace{-.3cm}

\noindent
 Key words and phrases: Polynomials, zeros of polynomials, Polar
Derivative.
 }
\newpage

It is noted that in (1.1) equality hold if and only if $p(z)$ has
all  its zeros at the origin and so it is natural to seek
improvements under appropriate assumptions on the zeros of $p(z)$.

If $p(z)$ having no zeros in $|z|<1$, then the above inequality can
be replaced by
\begin{equation}
\max_{|z|=1} |p'(z)|\leq \frac{n}{2} \max_{|z|=1} |p(z)|.
\end{equation}
Inequality (2) was conjectured by Erdos and later proved by Lax
\cite{L}. On the other hand, it was shown by Turan \cite{T} that if
all the zeros of $p(z)$ lie in $|z|< 1$, then
\begin{equation}
\max_{|z|=1} |p'(z)|\geq \frac{n}{2} \max_{|z|=1} |p(z)|.
\end{equation}
The above inequality was generalized by Govil \cite{G2}. Who proved
that if $p(z)$ is a polynomial of degree $n$ having all its zeros in
$|z|< k$, then for $k\leq 1$
\begin{equation}
\max_{|z|=1} |p'(z)|\geq \frac{n}{1+k} \max_{|z|=1} |p(z)|,
\end{equation}
and for $k\geq 1$
\begin{equation}
\max_{|z|=1} |p'(z)|\geq \frac{n}{1+k^n} \max_{|z|=1} |p(z)|.
\end{equation}
Both the above inequalities are best possible with equality in (4)
holding for $p(z) = (z+k)^n$, while in (5)the equality holds for the
polynomial $p(z) = z^n +k^n$. As an extension of (2) was shown by
Malik \cite{M} that, if $p(z)\ne 0$ in $|z|< k, k\geq 1$, then
\begin{equation}
\max_{|z|=1} |p'(z)|\leq \frac{n}{1+k} \max_{|z|=1} |p(z)|.
\end{equation}
Equality in (6) holds for $p(z) = (z+k)^n$.

By considering a more general class of polynomials $p(z) = a_0 +
\sum_{\nu =\mu}^na_\nu z^\nu$, $1\leq\mu\leq n,$  not vanishing in
$|z|< k$, $k>0$,  then for $0<r\leq R\leq k$, inequality (6) is
generalized by Aziz and Shah \cite{A4} by proving
\begin{equation}
\max_{|z|=R} |p'(z)|\leq nR^{\mu
-1}\frac{(R^{\mu}+k^{\mu})^{\frac{n}{\mu}-1}}{(r^{\mu}+k^{\mu})^{\frac{n}{\mu}}}
[ \max_{|z|=r} |p(z)| -\min_{|z|=k}|p(z)|].
\end{equation}
Equality in (7) holds for $p(z) = (z^\mu+k^\mu)^{\frac{n}{\mu}}$
where $n$ is a multiple of $\mu$. On the other hand, for the class
of polynomial $p(z) = a_nz^n + \sum_{\nu =
\mu}^{n}a_{n-\nu}z^{n-\nu}, 1\leq \mu\leq n,$ of degree $n$ having
all its zeros in $|z|\leq k, k\leq 1$, Aziz and Shah \cite{A2}
proved
\begin{equation}
\max_{|z|=1} |p'(z)|\geq \frac{n}{1+k^\mu} \{\max_{|z|=1}
|p(z)|+\frac{1}{k^{n-\mu}}\min_{|z|=k} |p(z)|\}.
\end{equation}
Let $D_\alpha \{p(z)\}$ denote the polar derivative of the
polynomial $p(z)$ of degree $n$ with respect to the point $\alpha$,
then $$ D_\alpha \{p(z)\} =np(z) + (\alpha -z)p'(z).$$ The
polynomial $D_\alpha \{p(z)\}$ is of degree at most $n-1$ and it
generalized the ordinary derivative in the sense that $$ \lim_
{\alpha\rightarrow\infty}[\frac{D_\alpha \{p(z)\}}{\alpha}] =p'(z)$$
 Dewan \cite{D2}
extended (8) to the polar derivative of a polynomial and proved the
following.
\paragraph{\bf Theorem A} Let
 $p(z) = a_nz^n + \sum_{\nu = \mu}^{n}a_{n-\nu}z^{n-\nu},
1\leq \mu\leq n,$ be a polynomial of degree $n$ having all its zeros
in $|z|\leq k, k\leq 1$. Then for every real or complex number
$\alpha$ with $|\alpha| \geq k^\mu$, we have
 \begin{align}
\hspace{0.4cm} \max_{|z|=1}|D_\alpha p(z)|&\geq \frac{n(|\alpha
|-k^\mu)}{1+k^\mu}
\max_{|z|=1} |p(z)|+n(\frac{|\alpha |+1}{k^{n-\mu}(1+k^\mu)})m\notag\\
&\hspace{0.5cm}n (\frac{k^\mu-A\mu}{1+k^\mu})\max_{|z|=1}|p(z)|
+\frac{n(A_\mu-k^\mu)}{k^n(1+k^\mu)}m,
\end{align}
where $m=\min_{|z|=k}|p(z)|$ and
$$A_\mu=\frac{n(|a_n|-\frac{m}{k^n})k^{2\mu}+\mu|a_{n-\mu}|k^{\mu-1}}{n(|a_n|-\frac{m}{k^n})k^{\mu-1}+\mu|a_{n-\mu}|}.$$
Dividing both sides of inequality (9) by $|\alpha |$ and let
$|\alpha |\rightarrow \infty$ we get (8).

As an extension of (6) to the polar derivative of a polynomial, we
have the following result due to Dewan \cite{D1}.
\paragraph{\bf Theorem B}
If  $p(z) = a_nz^n + \sum_{\nu = \mu}^{n}a_{n-\nu}z^{n-\nu}, 1\leq
\mu\leq n,$  is a polynomial of degree $n$ having no zeros in
$|z|<k, k\geq 1$, then for every real or complex number $\alpha$
with $|\alpha| \geq 1,$
\begin{equation}
\max_{|z|=1} |D_\alpha p(z)|\leq \frac{n}{1+k^\mu}\{(|\alpha
|+k^\mu) \max_{|z|=1} |p(z)|-(|\alpha |-1)m\}
\end{equation}
where $m=\displaystyle\min_{|z|=k}|p(z)|$ .\\ The above theorem is
an extension of a result of Aziz \cite{A1} and for $\mu=1$, it
reduces to a result of Aziz and Shah \cite{A2}.

In this paper, we first obtain the following generalization of
inequality (8) which is also a refinement of inequality (9).
\paragraph{\bf Theorem 1.}
Let
 $p(z) = a_nz^n + \sum_{\nu = \mu}^{n}a_{n-\nu}z^{n-\nu},
1\leq \mu\leq n,$ be a polynomial of degree $n$ having all its zeros
in $|z|\leq k, k\leq 1$ and $\alpha$ is any real or complex number
with $|\alpha| \geq \frac{k^\mu}{R^{\mu-1}},$ then for $rR\geq k^2$
and $r\leq R$, we have
\begin{align}
 %\hspace{0.4cm}
 \max_{|z|=R} |D_\alpha p(z)|&\geq n(R^{\mu-1}|\alpha
|-k^\mu)\frac{(R^\mu+k^\mu )^{\frac{n}{\mu
}-1}}{(r^{\mu}+k^{\mu})^{\frac{n}{\mu}}} \max_{|z|=r}
|p(z)|\notag\\
&+\frac{n(R^{\mu-1}|\alpha |+R^\mu)}{k^{n-\mu}(R^\mu+k^\mu)}
\min_{|z|=k}
|p(z)|\notag\\
&+n(k^\mu-R^\mu
A'_\mu)\frac{(R^{\mu}+k^{\mu})^{\frac{n}{\mu}-1}}{(r^{\mu}+k^{\mu})^{\frac{n}{\mu}}}
 \max_{|z|=r} |p(z)|\notag\\
&+\frac{nR^n(R^\mu A'_\mu-k^\mu)}{k^n(R^\mu+k^\mu)}\min_{|z|=k} |p(z)|\notag\\
&+\frac{nR^{\mu-1}}{(R^\mu+k^\mu)}(|\alpha|-RA'_\mu)[(\frac{R}{r})^n-(\frac{R^{\mu}+k^{\mu}}{r^{\mu}+k^{\mu}})^{\frac{n}{\mu}}]
 \min_{|z|=k} |p(z)|,
\end{align}
where
$$A'_\mu=\frac{n(|a_n|-\frac{m}{k^n})\frac{k^{2\mu}}{R^\mu}+\mu|a_{n-\mu}|\frac{k^{\mu-1}}{R^{\mu-1}}}{nR(|a_n|-\frac{m}{k^n})k^{\mu-1}+\mu|a_{n-\mu}|}$$
and $m=\displaystyle\min_{|z|=k}|p(z)|.$
\paragraph{\bf Remark 1.} For $R=r=1$ theorem 1 reduces to (9).
\paragraph{\bf Remark 2.}
 Dividing the two sides of (11) by $|\alpha |$, letting $|\alpha |\rightarrow \infty$, we
 get
\begin{align}
 %\hspace{0.4cm}
 \max_{|z|=R} |p'(z)|&\geq \frac{nR^{\mu-1}(R^\mu+k^\mu )^{\frac{n}{\mu
}-1}}{(r^{\mu}+k^{\mu})^{\frac{n}{\mu}}} \max_{|z|=r} |p(z)|
+\frac{nR^{\mu-1}}{k^{n-\mu}(R^\mu+k^\mu)}\min_{|z|=k} |p(z)|\notag\\
&+\frac{nR^{\mu-1}}{(R^\mu+k^\mu)}[(\frac{R}{r})^n-\frac{(R^{\mu}+k^{\mu})}{(r^{\mu}+k^{\mu})}^{\frac{n}{\mu}}]
\min_{|z|=k} |p(z)|.
\end{align}
This includes inequality (8) as special case.
%\bigskip
The following result immediately follows by taking $k=1$ in theorem
1.
\paragraph{Corollary 1.} If  $p(z) = a_nz^n + \sum_{\nu = \mu}^{n}a_{n-\nu}z^{n-\nu},
1\leq \mu\leq n,$ be a polynomial of degree $n$ having all its zeros
in $|z|\leq 1$ and $\alpha$   is any real or
 complex number $\alpha$ with $|\alpha| \geq \frac{1}{
R^{\mu-1}}$, then for $rR\geq 1$ and $r\leq R$, we have
\begin{align}
 %\hspace{0.4cm}
 \max_{|z|=R} |D_\alpha p(z)|&\geq \frac{n(R^{\mu-1}|\alpha
|+R^\mu)}{R^\mu+1} \min_{|z|=1}
|p(z)|\notag\\
&+n(R^{\mu-1}|\alpha |-1)\frac{(R^\mu+1 )^{\frac{n}{\mu}
-1}}{(r^{\mu}+1)^{\frac{n}{\mu}}} \max_{|z|=r}
|p(z)|\notag\\
%&+\frac{n(R^{m-1}|\alpha |-R^n)}{k^{n-m}(R^n+k^m)} \min_{|z|=k}
%|P(z)|\notag\\+\frac{n(R^{m-1}|\alpha |-R^n)}{k^{n-m}(R^n+k^m)}
%\min_{|z|=k}
%|P(z)|\notag\\
&+n(1-R^\mu
A''_\mu)\frac{(R^{\mu}+1)^{\frac{n}{\mu}-1}}{(r^{\mu}+1)^{\frac{n}{\mu}}}
 \max_{|z|=r} |p(z)|\notag\\
&+\frac{nR^n(R^\mu A''_\mu-1)}{(R^\mu+1)}\min_{|z|=1} |p(z)|\notag\\
&+\frac{nR^{\mu-1}}{(R^\mu+1)}(|\alpha|-RA''_\mu)[(\frac{R}{r})^n-\frac{(R^{\mu}+1)}{(r^{\mu}+1)}^{\frac{n}{\mu}}]
 \min_{|z|=1} |p(z)|,
\end{align}
where
$$A^{\prime\prime}_\mu=\frac{\frac{n}{R^\mu}(|a_n|-m)+\frac{\mu}{R^{\mu-1}}|a_{n-\mu}|}{nR(|a_n|-m)+\mu|a_{n-\mu}|}$$
and $m=\displaystyle\min_{|z|=1}|p(z)|.$\\ We next prove the
following result which is a generalization of the inequality (7),and
result due to Gardner, Govil and Weems \cite{G1}.
\paragraph{Theorem 2.} If  $p(z) = a_0 + \sum_{\nu = m}^{n}a_{\nu}z^{\nu},
1\leq m\leq n,$ be a polynomial of degree $n$ having no zeros in
$|z|< k, k\geq 1$ for $0\leq r\leq R\leq k$, then for every real or
complex number $\alpha$ with $|\alpha| \geq R,$
\begin{align}
 %\hspace{0.4cm}
 \max_{|z|=R} |D_\alpha p(z)|\leq &\notag\frac{n(R^{m-1}|\alpha
|+k^m)}{R^m+k^m}[(\frac{R^m+k^m}{r^m+k^m})^\frac{n}{m}\max_{|z|=r}
|p(z)|\notag\\
 &-((\frac{R^m+k^m}{r^m+k^m})^\frac{n}{m}-1)
  \min_{|z|=k}
|p(z)|]\notag\\
&-\frac{n(R^{m-1}|\alpha |-R^m)}{R^m+k^m}\min_{|z|=k} |p(z)|.
\end{align}
Equality in (14) holds for $p(z) = (z^m+k^m)^{\frac{n}{m}}$.
\paragraph{Remark 3.} For $R=r=1$ theorem 2 reduces to theorem B.

\paragraph{Remark 4.} Dividing the two sides of (14) by $|\alpha|$ and letting
$|\alpha |\rightarrow \infty$, we have the following inequality,
which is an improvement as well as a generalization of a result
proved by Bidkham and Dewan \cite{B}
\[
 \max_{|z|=R} |p'(z)|\leq \frac{nR^{m-1}(R^m+k^m)^{\frac{n}{m
}-1}}{(r^{m}+k^{m})^{\frac{n}{m}}} \{\max_{|z|=r}
|p(z)|-\min_{|z|=k} |p(z)|\}.\]
 The result is best possible and
equality holds for the polynomial $p(z) = (z^m+k^m)^{\frac{n}{m}}$
where $n$ is a multiple of $m$.

\section*{\small 2. LEMMAS}
For the proofs of these theorems we needs the following lemmas.
\paragraph{Lemma 2.1.}  If $p(z) = a_0 + \sum_{\nu = m}^{n}a_{\nu}z^{\nu},
1\leq m\leq n,$ is a polynomial of degree $n$ such that $p(z)\ne 0$
in  $|z|<k, k>0,$ then for $0\leq r\leq R\leq k,$
\begin{equation}
 \max_{|z|=r} |p(z)|\geq
 (\frac{r^m+k^m}{R^m+k^m})^\frac{n}{m}\max_{|z|=R}|p(z)|
+[1-(\frac{r^m+k^m}{R^m+k^m})^\frac{n}{m}
  ]\min_{|z|=k}
|p(z)|.
\end{equation}
 Here the result is best possible and equality hold for
the polynomial $p(z) = (z^m+k^m)^{\frac{n}{m}}$ where $n$  is
multiple of $m$. This result is due to Dewan [7].
\paragraph{Lemma 2.2}. If $p(z) = a_0 + \sum_{\nu = \mu}^{n}a_{\nu}z^{\nu},
1\leq \mu\leq n,$ is a polynomial of degree $n$ having no zeros in
$|z|\leq k, k\geq 1,$ then for every real or complex number $\alpha$
with $|\alpha |\geq 1$
\[ \max_{|z|=1} |D_\alpha
p(z)|\leq\frac{n}{1+k^\mu}\{(|\alpha|+k^\mu)\max_{|z|=1}
|p(z)|-(|\alpha|-1)m\},
\]
where $m=\displaystyle\min_{|z|=k}|p(z)|.$\\This result is due to
[6].
\section*{\small 3. PROOF OF THE THEOREMS}
{\bf Proof of theorem 1.} By hypothesis the polynomial $p(z) =
a_nz^n + \sum_{\nu = \mu}^{n}a_{n-\nu}z^{n-\nu}, 1\leq \mu\leq n,$
has all its zeros in $|z|\leq k,$ where $k\leq 1$, therefore it
follows that $F(z)=P(Rz)$ has all its zeros in $|z|\leq
\frac{k}{R},$ where $\frac{k}{R}\leq 1$ and $|\alpha |\geq
(\frac{k}{R})^\mu$. Applying inequality (9) to the polynomial $F(z)$
we get
\begin{align}
 %\hspace{0.4cm}
 \max_{|z|=1} |D_\alpha F(z)|&\geq \frac{n(|\alpha
|-\frac{k^\mu}{R^\mu})}{(1+\frac{k^\mu}{R^\mu})} \max_{|z|=1}
|F(z)|\notag\\
&+n(\frac{|\alpha|+1}{(\frac{k}{R})^{n-\mu}(1+\frac{k^\mu}{R^\mu})})
\min_{|z|=k/R}|F(z)|\notag\\
&+\frac{n((\frac{k}{R})^\mu-A'_\mu)}{(1+(\frac{k}{R})^\mu)}
\max_{|z|=1} |F(z)|
+\frac{n(A'_\mu-\frac{k^\mu}{R^\mu})}{(\frac{k}{R})^n(1+(\frac{k}{R})^\mu)}\min_{|z|=k/R}
|F(z)|,
\end{align}
where
$$A'_\mu=\frac{n(|a_n|-\frac{m}{k^n})\frac{k^{2\mu}}{R^\mu}+\mu|a_{n-\mu}|\frac{k^{\mu-1}}{R^{\mu-1}}}{nR(|a_n|-\frac{m}{k^n})
k^{\mu-1}+\mu|a_{n-\mu}|}$$ and
$m=\displaystyle\min_{|z|=k}|p(z)|,$\\which is equivalent to
\begin{align*}
& \max_{|z|=1}|np(Rz)+(\alpha-z)Rp'(Rz)|\geq\frac{n(|\alpha |R^\mu-k^\mu)}{R^\mu+k^\mu}\max_{|z|=1}|p(Rz)|\notag\\
&+\frac{nR^n(|\alpha |+1)}{K^{n-\mu}(R^\mu+K^\mu)} \min_{|z|=k/R}
|p(Rz)| +n\frac{(K^\mu-R^\mu A'_\mu)}{(R^\mu+k^\mu)} \max_{|z|=1} |p(Rz)|&\notag\\
&+\frac{nR^n(R^\mu A'_\mu-k^\mu)}{k^n(R^\mu+k^\mu)}\min_{|z|=k/R}
|p(Rz)|,
\end{align*}
which gives
\begin{align}
& \max_{|z|=R}|np(z)+(\alpha R-z)p'(z)|\geq\frac{n(|\alpha R|R^{\mu-1}-k^\mu)}{R^\mu+k^\mu}\max_{|z|=R}|p(z)|\notag\\
&+\frac{n(|\alpha R|R^{n-1}+R^n)}{K^{n-\mu}(R^\mu+K^\mu)}
\min_{|z|=k}
|p(z)| +n(\frac{K^\mu-R^\mu A'_\mu}{R^\mu+k^\mu}) \max_{|z|=R} |p(z)|&\notag\\
&+\frac{nR^n(R^\mu A'_\mu-k^\mu)}{k^n(R^\mu+k^\mu)}\min_{|z|=k}
|p(z)|.
\end{align}
On other hand,since $p(z)$ has all its zeros in $|z|<k, k>0$
 therefore it follows that $q(z) = z^n\overline{p(\frac{1}{z})}\ne 0$
 for $|z|<1/k.$
Applying inequality (15) to $q(z)$, we get
\[
 \max_{|z|=1/R} |q(z)|\geq
 (\frac{(\frac{1}{k})^\mu+(\frac{1}{R})^\mu}{(\frac{1}{k})^\mu+(\frac{1}{r})^\mu})^\frac{n}{\mu}\max_{|z|=1/r}|q(z)|
+[1-(\frac{(\frac{1}{k})^\mu+(\frac{1}{R})^\mu}{(\frac{1}{k})^\mu+(\frac{1}{r})^\mu})^\frac{n}{\mu}
  ]\min_{|z|=1/k}
|q(z)|
\],
which gives
\begin{equation}
 \max_{|z|=R} |p(z)|\geq
 (\frac{R^\mu+k^\mu}{r^\mu+k^\mu})^\frac{n}{\mu}\max_{|z|=r}|p(z)|
+[(\frac{R}{r})^n-(\frac{R^\mu+k^\mu}{r^\mu+k^\mu})^\frac{n}{\mu}
  ]\min_{|z|=k}
|p(z)|.\end{equation} From (17) and (18), we have
\begin{align*}
 %\hspace{0.4cm}
 \max_{|z|=R} |D_{R\alpha} p(z)|&\geq n(|\alpha
 R|R^{\mu-1}-k^\mu)\frac{(R^\mu+k^\mu )^{\frac{n}{\mu
}-1}}{(r^{\mu}+k^{\mu})^{\frac{n}{\mu}}} \max_{|z|=r}
|p(z)|\notag\\
&+\frac{n(R^{n-1}|\alpha R|+R^n)}{k^{n-\mu}(R^\mu+k^\mu)}
\min_{|z|=k}
|p(z)|\notag\\
&+n(k^\mu-R^\mu
A'_\mu)\frac{(R^{\mu}+k^{\mu})^{\frac{n}{\mu}-1}}{(r^{\mu}+k^{\mu})^{\frac{n}{\mu}}}
 \max_{|z|=r} |p(z)|\notag\\
&+\frac{nR^n(R^\mu A'_\mu-k^\mu)}{k^n(R^\mu+k^\mu)}\min_{|z|=k} |p(z)|\notag\\
&+\frac{nR^{\mu-1}}{(R^\mu+k^\mu)}(|R\alpha|-RA'_\mu)[(\frac{R}{r})^n-(\frac{R^{\mu}+k^{\mu}}{r^{\mu}+k^{\mu}})^{\frac{n}{\mu}}]
 \min_{|z|=R} |p(z)|,
\end{align*}
which is equivalent to (11).

\bigskip

\noindent{\bf Proof of theorem 2.} By hypothesis the polynomial
$p(z) = a_0+ \sum_{\nu = m}^{n}a_{\nu}z^{\nu}, 1\leq m\leq n,$
having no zeros in $|z|< k,$ where $k\geq 1$, therefore it follows
that $F(z)=p(Rz)$ having no zeros in $|z|\leq \frac{k}{R},$ where
$\frac{k}{R}\geq 1$. Hence using lemma 2 for $|\alpha|\geq1$. we
have
\[ \max_{|z|=1} |D_\alpha
F(z)|\leq\frac{n}{1+(\frac{k}{R})^m}\{(|\alpha|+(\frac{k}{R})^m)\max_{|z|=1}
|F(z)|-(|\alpha|-1)\min_{|z|=k/R} |F(z)|\},
\]
which is equivalent to
\[\max_{|z|=1}|np(Rz)+(\alpha-z)Rp'(Rz)|\leq \frac{n}{1+(\frac{k}{R})^m}\{(|\alpha|+(\frac{k}{R})^m)\max_{|z|=1}
|p(Rz)|-(|\alpha|-1)\min_{|z|=k/R} |p(Rz)|\}.
\]
Replacing $Rz$ by $z$ we get
\begin{align*}
 %\hspace{0.4cm}
 \max_{|z|=R} |D_{R\alpha}p(z)|\leq &\notag\frac{n(R^{m-1}|\alpha R
|+k^m)}{R^m+k^m}\max_{|z|=R}
|p(z)|\notag\\
 &-
\frac{n(R^{m-1}|R\alpha |-R^m)}{(R^m+k^m )}\min_{|z|=k} |p(z)|
\end{align*}
for $R\geq 1$\\ or \begin{align*}
 %\hspace{0.4cm}
 \max_{|z|=R} |D_\alpha p(z)|\leq &\notag\frac{n(R^{m-1}|\alpha
|+k^m)}{R^m+k^m}\max_{|z|=R}
|p(z)|\notag\\
 &-
\frac{n(R^{m-1}|\alpha |-R^m)}{(R^m+k^m)}\min_{|z|=k}
|p(z)|.~~~~~~~~~~~~~~~~~~~~~~~~~~~~~~~~~~~~~~~~~~~~~(19)
\end{align*}

 By combining (19) and (15) we get (11).

%%%%%%%%%%%%%%reference


\begin{thebibliography}{99}






\bibitem{A1} A. AZIZ, Inequalities for the polar derivative  of a
polynomial, J. Approx. Theory 55 (1988) 183--193.

\bibitem{A2} A. AZIZ AND W.M. SHAH, An integral mean estimate for
polynomial, Indian J. Pure Apple. Math. 28 (1997) 1413--1419.

\bibitem{A3} A. AZIZ AND W.M. SHAH, Inequalities for polar derivative of a polynomial,
Indian J. Pure Appl. Math. 29 (1998) 163--173.

\bibitem{A4} A. AZIZ AND W.M. SHAH, Inequalities for a polynomial and its
derivative, Math. Inequal. Appl.7(3)(2004) 379--391.

\bibitem{B} M. BIDKHAM AND K.K. DEWAN, Inequalities for polynomial and
its derivative, J. Math. Anal. Appl. 166 (1992)319--324.

\bibitem{D1} K.K. DEWAN AND NARESH SINGH and ABDULLAH MIR, Extensions of some polynomial,
inequalities to the polar derivative,  J. Math. Anal. Appl. 352
(2009), no. 2, 807--815.

\bibitem{D2} K.K. DEWAN, R.S. YADAV AND M.S.PUKHTA, Inequalities for a polynomial
and its derivative, Math. Inequal. Appl. 2(2)  (1999)203--205.

\bibitem{G1}R.B. GARDNER, N.K. GOVIL AND A.WEEMS, Some result concerning
rate of growth of polynomials, East J. Apporox. 10(2004) 301--312.

\bibitem{G2} N.K. GOVIL, On the derivative of a polynomial, Proc. Amer. Math. Soc. 41(1973) 643--546.

\bibitem{L} P.D. LAX, Proof of a conjecture of P.Erdos on the derivative of
a polynomial, Bull. Amer. Math. Soc. 50 (1944) 509--513.

\bibitem{M} M.A. MALIK, On the derivative of a
polynomial, J. Londan Math. Soc. 1(1969) 57--60.

\bibitem{T} P. TURAN, Uberdie ableting von polynoment, Compos. Math. 7(1939)
82--95.



\end{thebibliography}
\end{document}